\documentclass[12pt]{amsart}
\usepackage{etex}
\makeindex

\usepackage{graphicx,tikz,url,hyperref}
\usepackage[latin1]{inputenc} \usepackage[T1]{fontenc} \usepackage{lmodern}
\usepackage{pgf}
\usepackage{epsfig,graphicx,color,amsmath,a4wide,amssymb,amsfonts,amsbsy,xy,latexsym,epsf,pstricks,verbatim,times}
\usepackage{MnSymbol}
\usepackage{pgf,color}

\usepackage{changebar}

\definecolor{LightGrey}{rgb}{.85,.85,.85}
\definecolor{DarkGrey}{rgb}{.5,.5,.5}

\definecolor{Blue}{rgb}{.0,.0,0.9}
\definecolor{LightBlue1}{rgb}{.2,.4,0.9}
\definecolor{LightBlue2}{rgb}{.3,.5,0.9}
\definecolor{LightBlue3}{rgb}{.4,.6,0.9}
\definecolor{LightBlue4}{rgb}{.5,.7,.9}
\definecolor{LightBlue5}{rgb}{.6,.8,.9}
\definecolor{LightBlue6}{rgb}{.7,.9,.9}

\definecolor{Red}{rgb}{.9,.0,.0}
\definecolor{LightRed1}{rgb}{0.9,.2,.4}
\definecolor{LightRed2}{rgb}{0.9,.3,.5}
\definecolor{LightRed3}{rgb}{0.9,.4,.6}
\definecolor{LightRed4}{rgb}{.9,.5,.7}
\definecolor{LightRed5}{rgb}{.9,.6,.8}
\definecolor{LightRed6}{rgb}{.9,.7,.9}

\def\KK{{\mathbf K}}

\definecolor{Grey}{rgb}{.5,.5,.5}

\definecolor{Blue}{rgb}{.0,.0,0.9}
\definecolor{LightBlue1}{rgb}{.2,.4,0.9}
\definecolor{LightBlue2}{rgb}{.3,.5,0.9}
\definecolor{LightBlue3}{rgb}{.4,.6,0.9}
\definecolor{LightBlue4}{rgb}{.5,.7,.9}
\definecolor{LightBlue5}{rgb}{.6,.8,.9}
\definecolor{LightBlue6}{rgb}{.7,.9,.9}

\definecolor{Red}{rgb}{.9,.0,.0}
\definecolor{LightRed1}{rgb}{0.9,.2,.4}
\definecolor{LightRed2}{rgb}{0.9,.3,.5}
\definecolor{LightRed3}{rgb}{0.9,.4,.6}
\definecolor{LightRed4}{rgb}{.9,.5,.7}
\definecolor{LightRed5}{rgb}{.9,.6,.8}
\definecolor{LightRed6}{rgb}{.9,.7,.9}

\xyoption{all}
\usepackage[english]{babel}

\newcounter{noalgo}[section]
\newdimen\indentalgo
\newdimen\indentalgodec\indentalgo=0.0mm\indentalgodec=10mm

\newcommand{\If}{\advance\indentalgo by \indentalgodec {\bf if }}

\newcommand{\For}{\global\advance\indentalgo by \indentalgodec {\bf for }}

\newcommand{\Endindent}{\global\advance\indentalgo by -\indentalgodec}

\newdimen\decalage \decalage=0.5cm
\newcounter{algo} 

\newcommand{\PP}{\mathbf P}

\def\<<{\leavevmode
  \raise0.28ex\hbox{$\scriptscriptstyle\langle\!\langle$}\nobreak
  \hskip -.6pt plus.3pt minus.2pt\,}
\def\>>{\,\nobreak\hskip -.6pt plus.3pt minus.2pt
  \raise0.28ex\hbox{$\scriptscriptstyle\rangle\!\rangle$}}

\def\AA{{\mathbf A}}

\def\vol{{{v}}}

\def\Proj{\mathop{\rm{Proj}}\nolimits }
\def\Hom{\mathop{\rm{Hom}}\nolimits }

\def\Ker{\mathop{\rm{Ker}}\nolimits }

\def\Spec{\mathop{\rm{Spec}}\nolimits }

\def\CC{{\mathbf C}}
\def\OO{{\mathbf {O}}}

\def\AA{{\mathbf A}}

\def\QQ{{\mathbf Q}}

\def\RR{{\mathbf R}}

\def\ZZ{{\mathbf Z}}
\def\Dgot{{\mathfrak D}}

\def\cI{{\mathcal I}}

\def\cL{{\mathcal L}}

\def\cM{{\mathcal M}}

\def\cO{{\mathcal O}}

\newtheorem{proposition}{Proposition}
\newtheorem{theorem}{Theorem}

\providecommand{\myproofname}{Proof}

\usepackage{enumitem}

\begin{document}

\begin{abstract}
We construct small models of number fields and deduce a better bound for the number of number fields of given degree and  bounded discriminant.
\end{abstract}

\title{Enumerating number fields}

\author{Jean-Marc Couveignes}
\address{Jean-Marc Couveignes, Univ. Bordeaux, CNRS,
  Bordeaux-INP, IMB, UMR 5251, F-33400 Talence, France.}
\address{Jean-Marc Couveignes, INRIA, F-33400 Talence, France.}
\email{Jean-Marc.Couveignes@u-bordeaux.fr}

\date{\today}

\maketitle
\setcounter{tocdepth}{2} 
\tableofcontents

\section{Introduction}

We prove the two theorems below.

\begin{theorem}[Number fields have small models]\label{th:sm}
  There exists a positive constant $\cO$ such that the following is true.
  Let $\KK$ be a number field of degree $n\geqslant \cO$ and discriminant $d_\KK$ over $\QQ$.
  There exist  integers $r\leqslant \cO \log n$
  and $d\leqslant \cO \log n$
such that ${d+r\choose r}\leqslant \cO n\log n$ and there exists 
    $r$ polynomials $E_1$, $E_2$,
  \dots , $E_r$
  of degree $\leqslant d$ in $\ZZ[x_1, \ldots, x_r]$ all having  coefficients bounded in absolute
  value by $n^{\cO\log n}d_\KK^{\cO\frac{\log n}{n}}$ such that 
    the (smooth and zero-dimensional affine) scheme with equations
    \[E_1 = E_2 = \dots = E_r=0 \text{ and } \det \left( \partial E_i/\partial x_j \right)_{1\leqslant i, \, j\leqslant r} \not = 0\]
    contains $\Spec \KK$ as one of its irreducible components.
\end{theorem}

\begin{theorem}[Number fields with bounded discriminant]\label{th:nf}
  There exists a positive constant $\cO$ such that the following is true.
  Let $n\geqslant \cO$ be an integer. Let $H\geqslant 1$ be an integer.
  The number of isomorphism classes of number fields with degree $n$ and discriminant
  $\leqslant H$
  is
  $\leqslant n^{\cO n\log^3n}H^{\cO \log^3 n}$.
\end{theorem}

The meaning of Theorem \ref{th:sm} is that we can describe a number field using few parameters in some sense.
We have a short description of it as a  quotient  of a finite  algebra : the smooth zero-dimensional
part of a complete intersection of small degree and small height in a projective space of small dimension.

Theorem \ref{th:nf} improves on previous results by Schmidt \cite{Sch} and Ellenberg-Venkatesh \cite{EV}.
Schmidt obtains a bound $H^{\frac{n+2}{4}}$ times a function of  $n$.
 Ellenberg and Venkatesh obtain a bound $H^{\exp (\cO \sqrt{\log n})}$ times a function of $n$. 
We combine techniques from geometry of numbers and   interpolation  theory  to produces small
projective models of $\Spec \KK$ and lower  the exponent of $H$ down to $\cO \log^3 n$. A key point is to
look for local equations rather than a full set of generators of the ideal of these models.

Our estimate is not sharp of course. Indeed for $n=1$
the exponent of $H$
can be taken to be  $0$.
For $2\leqslant n\leqslant 5$ the exponent of $H$ can be taken to be $1$ according to work by Davenport and Heilbronn \cite{DH}
for $n=3$, and Bhargava \cite{Bh4, Bh5} for $n=4, \, 5$. It is a bit
delicate  to infer a general conjecture from
these results because
the techniques used for these small values of $n$ seem to be quite specific.
Cohen, Diaz and Olivier have collected experimental data
e.g. in \cite{Co1, Co2, CDO} suggesting that
 the number
of isomorphism classes of number fields of degree $n$
and discriminant $\leqslant H$ should grow linearly in $H$
for fixed $n\geqslant 2$.
Malle has stated in \cite{Malle}
a   more general and accurate conjecture  on the distribution of Galois groups
of number fields that  would confirm this intuition.

In Section \ref{sec:short} we recall notation,  definitions and  elementary results from the geometry of numbers.
In Section \ref{sec:small} we construct models for number fields as irreducible components of complete intersections with small
height in low dimensional projective spaces.
The last section is devoted to the proof of Theorem \ref{th:sm} and Theorem \ref{th:nf}.

The author thanks Pascal Autissier, Karim Belabas, Georges Gras
and Christian
Maire for their comments and suggestions.

\section{Short integers}\label{sec:short}

Let $\KK$ be a number field and let $n$ be the degree of
$\KK$ over $\QQ$. Let $\OO$ be the ring of integers of $\KK$.
Let $(\rho_i)_{1\leqslant i\leqslant r}$ be the $r$ real embeddings of $\KK$.
Let $(\sigma_j, \bar\sigma_j)_{1\leqslant j\leqslant s}$ be the $2s$
complex
embeddings of $\KK$. We also denote by $(\tau_k)_{1\leqslant k\leqslant n}$ the $n=r+2s$
embeddings of $\KK$.
Let \[\KK_\RR = \KK\otimes_\QQ\RR = \RR ^r \times \CC^s\]
be the Minkowski space. We follow the presentation in \cite[Chapitre 1, \S 5]{Neukirch}.
An element $x$ of $\KK_\RR$ can be given  by
$r$ real components $(x_\rho)_\rho$
and $s$ complex components $(x_\sigma)_\sigma$. So we write
$x = ((x_\rho)_\rho , (x_\sigma)_\sigma)$.
For such  an $x$ in $\KK_\RR$ we denote by $||x||$ the
maximum of the absolute values of its $r+s$ components. 
The canonical metric on $\KK_\RR$ is defined
by \[<x,y> = \sum_{1\leqslant i\leqslant r}x_iy_i+
\sum_{1\leqslant j\leqslant s}x_j\bar y_j+\bar x_jy_j.\]
In particular the contribution of complex embeddings is counted twice
\[<x,x> =  \sum_{1\leqslant i\leqslant r}x_i^2+
2\sum_{1\leqslant j\leqslant s}|x_j|^2.\]
The corresponding Haar measure is said to be canonical also.
The canonical measure  of the convex body $\{x,  ||x||\leqslant 1\}$ is
\[2^r(2\pi )^{s}\geqslant 2^n.\]
The map $a\mapsto a\otimes 1$  injects $\KK$ and $\OO$ into
$\KK_\RR$. For every non-zero $x$ in $\OO$ we have \[||x||\geqslant 1.\]
Let $(\alpha_i)_{1\leqslant i\leqslant n}$ be any $\ZZ$-basis of $\OO$.
Set $A  = (\tau_j (\alpha_i))_{1\leqslant i, j \leqslant n}$. The product $A\bar A^t$
is the Gram matrix $B = (<\alpha_i, \alpha_j>)_{1\leqslant i, j \leqslant n}$ of the canonical form
in the basis $(\alpha_i)_i$.
This is a real symmetric positive matrix.
The volume of $\OO$  according to the canonical Haar measure
is  \[\vol_\OO = \sqrt{\det (B)} = |\det (A)|.\]
The square of the volume of $\OO$  is the discriminant of $\KK$
  \[d_\KK = \det (B) = |\det (A)|^2=\vol_\OO^2.\]
Applying Minkowski's second theorem \cite[Lecture III, \S 4, Theorem 16]{Siegel}
to the gauge function $x\mapsto ||x||$ we find that $\OO$ contains $n$ linearly independant
elements $\omega_1$, $\omega_2$, \ldots, $\omega_n$ such that
\[\prod_{1\leqslant i\leqslant n}||\omega_i||\leqslant \vol_\OO=d_\KK^{1/2}.\]
We assume that the sequence $i\mapsto ||\omega_i||$
is non-decreasing and deduce that
\[||\omega_i||\leqslant v_\OO^{1/(n+1-i)}\]for every $1\leqslant i\leqslant n$.
This inequality is a bit unsatisfactory because it provides  little
information on the largest $\omega_i$. To improve on this estimate
we use the fact that $\OO$ is an integral domain.
We let $m = \lceil (n+1)/2\rceil$ be the smallest integer bigger than
$n/2$. On the one hand
\[||\omega_i||\leqslant d_\KK^{1/n}\]for every $1\leqslant i\leqslant m$.
On the other hand the products \[(\omega_i\omega_j)_{1\leqslant i, j\leqslant m}\]
generate a  $\ZZ$-module of rank $n$. Otherwise  there would exist
a non-zero linear form $f : \OO \rightarrow \ZZ$ vanishing on these products.
So the $m$ forms $f\circ \omega_i$ would be orthogonal to the
$m$ vectors $\omega_j$. Then $m+m\leqslant n$. A contradiction.
We deduce that all the successive minima of $\OO$ are \[\leqslant d_\KK^{\, 2/n}.\]
In other words $\OO$ is well balanced.
\begin{proposition}[Number fields have small integers]\label{prop:balanced}
  The ring of integers $\OO$ of a number field $\KK$ with degree $n$ and
  discriminant $d_\KK$ contains $n$ linearly independant elements
  $(\alpha_i)_{1\leqslant i\leqslant n}$ over $\ZZ$ such that all the absolute
  values of all the $\alpha_i$ are $\leqslant d_\KK^{\, 2/n}$.
\end{proposition}

Bhargava, Shankar,  Taniguchi,  Thorne,  Tsimerman, and  ZhaoSee
prove in \cite{Bha}[Theorem 3.1]  a similar statement
which is  somewhat  stronger
but less accurate. 

\section{Small  models}\label{sec:small}

Let \[\KK_\CC = \KK\otimes_\QQ\CC = \CC ^n.\]
Let $d\geqslant 5$ and $r\geqslant 1$ be two integers. We assume that
\[n(r+1)\leqslant {d+r\choose d}.\]
Let $M$ be the set of  monomials of total degree $\leqslant d$
in the $r$ variables $x_1$, \ldots, $x_r$.
We have  \[\AA^r_\CC = \Spec \CC[x_1, \ldots, x_r]
\subset \Proj \CC[x_0, x_1, \ldots, x_r] = \PP^r_\CC.\]
Let $V_\CC$ be the $\CC$-linear space generated by $M$.
We may associate to every element in $M$ the corresponding
degree $d$ monomial in the $r+1$ variables
$x_0$, $x_1$, \ldots, $x_r$.
We thus identify $V_\CC$ with $H^0(\cO_{\PP^r_\CC}(d))$, the space  of homogeneous polynomials of degree $d$.

Let $(P_\tau)_\tau$
be $n$ pairwise distinct  points in \[\CC^r = \AA^r(\CC).\]
The $P_\tau$ are indexed by the $n$ embeddings of $\KK$.
These $n$ points form a set (a reduced zero-dimensional subscheme
of $\PP^r_\CC$) called  $P$.
We call $\cI$ the corresponding ideal sheaf on $\PP^r_\CC$.
We denote by $2P$ the scheme  associated
with $\cI^2$. It consists of
$n$ double points.
We say 
that the scheme  $2P$ is well poised (or non-special) in degree $d$
if it imposes $n(r+1)$ independent conditions on degree $d$
homogeneous polynomials. Equivalently, the map
\[H^0(\cO_{\PP^r}(d))\rightarrow H^0(\cO_{2P}(d))\]
is surjective.
This is the case if and only if the
$n(r+1) \times {d+r\choose d}$
matrix 
\[\cM_P^1 = [(m(P_\tau))_{ \tau ,  \, m\in M  },
  (\partial m/\partial x_1 (P_\tau))_{ \tau ,\, m\in M},
  (\partial m/\partial x_2 (P_\tau))_{\tau , \, m\in M }, \dots, 
  (\partial m/\partial x_r (P_\tau))_{\tau , \, m\in M}]
\] has maximal rank $n(r+1)$. We note
that $\cM_P^1$ consists of $r+1$ blocks of size
$n\times {d+r\choose d}$ piled vertically.
    It has maximal rank  for a generic $P$ when  $d\geqslant 5$, according to a theorem
    of Alexander \cite{Alexander}, generalized by Alexander  and Hirschowitz \cite{AH}. Chandler  \cite[Theorem 1]{Chandler}
    provides  a simpler statement and proof. The recent
    exposition and simplification by Brambilla and Ottaviani
\cite{Bram}     is very useful also.

    We now let  $(\alpha_i)_{1\leqslant i\leqslant n}$ be $n$ linearly independant
    short elements in $\OO$ as in Proposition \ref{prop:balanced}.
We pick $rn$ rational integers $(u_{i,j})_{1\leqslant i\leqslant n, \, 1\leqslant j \leqslant r}$ and we
set
\[\kappa_j =  \sum_{1\leqslant i\leqslant n}u_{i,j}\alpha_{i}\]
for $1\leqslant j\leqslant r$.
Let \[\epsilon_\QQ  :  \QQ[x_1, \ldots, x_r] \rightarrow \KK\] be the homomorphism of $\QQ$-algebras  sending $x_j$ to $\kappa_j$
for $1\leqslant j \leqslant r$. Let \[e_\QQ  : \Spec \KK \rightarrow  \AA^r_\QQ
\subset \PP^r_\QQ\] be the corresponding morphism of schemes.
Tensoring $\epsilon_\QQ$ by $\RR$ we obtain
an homomorphism
\[\epsilon_\RR  :  \RR[x_1, \ldots, x_r] \rightarrow \KK_\RR\] sending
$x_j$ to $((\rho(\kappa_j))_\rho, (\sigma (\kappa_j))_\sigma)$.
We call  \[e_\RR  : \Spec \KK_\RR \rightarrow  \AA^r_\RR \subset \PP^r_\RR \]
the corresponding morphism of schemes.
We define 
\[\epsilon_\CC  :  \CC[x_1, \ldots, x_r] \rightarrow \KK_\CC\] and
\[e_\CC  : \Spec \KK_\CC \rightarrow  \AA^r_\CC \subset \PP^r_\CC \]
similarly. In particular $\epsilon_\CC$ maps $x_j$ onto $(\tau (\kappa_j))_\tau$.

We now consider  the points  $(P_\tau)_\tau$
such that $x_0(P_\tau)=1$ and \[(x_j(P_\tau))_\tau=  \left( \sum_{1\leqslant i\leqslant n}u_{i,j}\tau(\alpha_{i}) \right)_\tau ,\]
for $1\leqslant j\leqslant r$
or equivalently \[P_\tau = (        \sum_{1\leqslant i\leqslant n}u_{i,j}\tau(\alpha_{i})       )_{1\leqslant j\leqslant r}\in \CC^r = \AA^r(\CC)\subset \PP^r(\CC).\]
The maximal minors of the corresponding matrix $\cM_P^1$ are polynomials of total degree $\leqslant dn(r+1)$ in the $u_{i,j}$ and one of them is not identically zero.
The latter determinant  cannot vanish on the cartesian product $[0,dn(r+1)]^{nr}$. Thus there exist $nr$ rational  integers $u_{i,j}$ in the range
\[[0,dn(r+1)]\] such that the corresponding scheme $2P$ is well poised. 
We   assume that  the $u_{i,j}$ meet these conditions.

Since $2P$ is well poised, $P$ is well poised also. So $e_\QQ$, $e_\RR$ and $e_\CC$ are
 closed embeddings. In order to describe them  efficiently we
 look for polynomials with degree
 $\leqslant d$ and small integer coefficients vanishing at $P$.
 We denote by $V_\RR = \RR[x_1, \ldots, x_r]_d$  the $\RR$-vector space
 of polynomials in $\RR[x_1, \ldots, x_r]$
 of degree $\leqslant d$. There is a unique $\RR$-bilinear
 form on $V_\RR$ that turns the set $M$
 of monomials into an orthonormal basis. The lattice of relations
 with integer coefficients and degree $\leqslant d$ is the intersection
 between $\Ker \epsilon_\RR$ and \[V_\ZZ = \ZZ [x_1, \ldots, x_r]_d.\] This is
 a free $\ZZ$-module $\cL \subset V_\RR$
 of  rank
 \[\ell = {d+r\choose d} -n.\] We set $L = \cL\otimes_\QQ \RR$
 the underlying  $\RR$-vector space
 and $L^\perp$ its orthogonal complement in $V_\RR$. We denote by $\cL^\perp$
 the intersection $\cL^\perp  = L^\perp \cap V_\ZZ$.
 Since $V_\ZZ$ is unimodular,
 $\cL$ and $\cL^\perp$ have the same
 volume. See  \cite[Corollary 1.3.5.]{Martinet}.
 We denote by $\hat \OO = \Hom (\OO,\ZZ)$ the dual
 of $\OO$, the ring of integers of $\KK$,  as a $\ZZ$-module. We call \[\epsilon_{\ZZ, d} :
 \ZZ [x_1, \ldots, x_r]_d\rightarrow \OO\] the evaluation map
 in degree $\leqslant d$.
 We observe that $\cL^\perp$ contains
 the image of $\hat \OO$ by the transpose
 map \[\hat \epsilon_{\ZZ, d} : \hat \OO \rightarrow
 \ZZ [x_1, \ldots, x_r]_d\]
 \noindent where we have identified $\ZZ [x_1, \ldots, x_r]_d$
 with its dual thanks to the canonical bilinear form.
 So the volume of $\cL$ is bounded from above by
 the volume of  $\hat \epsilon_{\ZZ, d} (\hat \OO)$.
 We consider the matrix
 \[\cM_P^0 = [(m(P_\tau))_{\tau , \, m\in M }]\]
 of the map $\epsilon_{\CC, d} = \epsilon_{\ZZ, d}\otimes_\ZZ\CC$
 in the canonical bases.
 If we prefer to use an integral basis of $\OO$ on the right
 we should multiply $\cM_P^0$ on the left 
 by the inverse $T$ of the matrix  of a basis of $\OO$ in the
 canonical basis. We deduce that the square of  the volume
 of $\hat \epsilon_{\ZZ, d} (\hat \OO)$ is the determinant
 of $T\cM_P^0(\cM_P^0)^tT^t$.
 Since $T\cM_P^0$ has real coefficients we have
 \[\det (T\cM_P^0(\cM_P^0)^tT^t) = \det \left(T\cM_P^0\left( \overline{\cM_P^0}\right)^t\bar T^t\right)
 =\det \left(\cM_P^0   \left( \overline{\cM_P^0}\right)^t  \right) / d_\KK.\]
 So the square of the volume of the lattice of relations is bounded by
 the determinant of the hermitian positive definite matrix
 $\cM_P^0\left( \overline{\cM_P^0}\right)^t$ 
 divided by $d_\KK$.

 Recall that the coefficients in $\cM_P^0$
 are degree $\leqslant d$ monomials in the $\kappa_j =  \sum_{1\leqslant i\leqslant n}u_{i,j}\alpha_{i}$. The coefficients $u_{i,j}$ are bounded form above
 by $dn(r+1)$. All the absolute values of the $\alpha_i$ are bounded
 from above by $d_\KK^{2/n}$. So the coefficients in $\cM_P^0$
 are bounded from above by \[(n^2d(r+1))^dd_\KK^{2d/n}.\]
 The coefficients in $\cM_P^0\left( \overline{\cM_P^0}\right)^t$ are bounded from
 above by \[\Dgot = {d+r\choose d}(n^2d(r+1))^{2d}d_\KK^{4d/n}.\]
 The  matrix $\cM_P^0\left( \overline{\cM_P^0}\right)^t$
 being hermitian positive
 definite, its determinant is bounded
 from above by the product of the diagonal terms.
 We deduce that the volume of the lattice $\cL$
 of relations  is bounded from above by $\Dgot ^{n/2}$.
Recall that the dimension of $\cL$ is \[\ell = {d+r\choose d}-n.\]
For any $x$  in $V_\RR$ we denote by $||x||$ the $\ell_2$-norm in
the monomial basis. The volume of the sphere  $\{x \in L,  ||x||\leqslant 1\}$ is $\geqslant 2^{\ell}\ell^{-\ell/2}$.
Applying Minkowski's second theorem \cite[Lecture III, \S 4, Theorem 16]{Siegel}
to the gauge function $x\mapsto ||x||$ we find that $\cL$ contains $\ell$ linearly independant
elements $E_1$, $E_2$, \ldots, $E_\ell$ such that
\[\prod_{1\leqslant i\leqslant \ell}||E_i||\leqslant \ell^{\ell/2}\Dgot^{n/2}.\]
We assume that the sequence $i\mapsto ||E_i||$
is non-decreasing and deduce that the size of the $i$-th equation is bounded from above
\[||E_i||\leqslant \ell^{\frac{\ell}{2(\ell+1-i)}} \Dgot^{\frac{n}{2(\ell+1-i)}}\]for every $1\leqslant i\leqslant \ell$.
Again, this inequality is a bit unsatisfactory because it provides  little
information on the largest equations. This time we see no other way around than forgetting
the last $n-1$ equations.
On the one hand \[ ||E_i||\leqslant \ell^{\ell/2n} \Dgot^{1/2}\]
for every $1\leqslant i\leqslant \ell +1 -n$.

On the other hand the scheme $2P$ is well poised and the  $\CC$-vector space generated
by the $E_i$ for
$1\leqslant i\leqslant \ell+1-n$ has codimension $n-1 < n$ in $L\otimes_\RR\CC$. 
So there exists at least one embedding $\tau$ such that the $(\ell +1 -n)\times r$ matrix
\[\left((\partial E_i /  \partial x_j)(P_\tau)\right)_{1\leqslant i\leqslant \ell+1-n, \,  1\leqslant j\leqslant r}\]
has maximal rank $r$. In more geometric terms the $\CC$-vector space  generated
by the $\ell+1-n$ first equations $(E_i)_{1\leqslant i\leqslant \ell +1-n}$  surjects onto
the cotangent  space to $\PP^r_\CC$ at the geometric point $P_\tau$ for at least one $\tau$.
This means that there exist $r$ integers
$1\leqslant i_1 < i_2 < \dots < i_r \leqslant \ell +1-n$
such that the minor determinant 
\[\det \left(  (\partial E_{i_k} /  \partial x_j)(P_\tau)\right)_{1\leqslant k,\,  j \leqslant r}\]
is non-zero for some $\tau$ and thus for all $\tau$ by Galois action.  

\begin{proposition}[Number fields have small models]\label{prop:sm}
  Let $\KK$ be a number field of degree $n$ and discriminant $d_\KK$ over $\QQ$.
  Let $d\geqslant 5$ and $r\geqslant 1$ be rational  integers such that
  \[n(r+1)\leqslant {d+r\choose d}.\] There exists $r$ polynomials $E_1$, $E_2$,
  \dots , $E_r$
  of degree $\leqslant d$ in $\ZZ[x_1, \ldots, x_r]$ having  coefficients bounded in absolute
  value by \[\ell^{\ell/2n}   \times {d+r\choose d}^{1/2}(n^2d(r+1))^{d}d_\KK^{2d/n}\] where
\[\ell = {d+r\choose d}-n,\]  
  and such that
    the (smooth and zero-dimensional affine) scheme with equations
    \[E_1 = E_2 = \dots = E_r=0 \text{ and } \det \left( \partial E_i/\partial x_j \right)_{1\leqslant i, \, j\leqslant r} \not = 0\]
    contains $\Spec \KK$ as one of its irreducible components.
\end{proposition}

\section{Proof of main results}\label{sec:proof}

In this section, the
notation $\cO$ stands for a positive absolute constant.
Any sentence  containing
this symbol becomes true if the symbol is replaced in every occurrence by some 
large enough real number.

We specialize  the values of the parameters $r$ and $d$ in Proposition \ref{prop:sm}. We will take $d=r$. It is evident
that ${2r \choose r}\geqslant 2^r$ so
\[\frac{1}{r+1}{2r \choose r}\geqslant 2^{\frac{r}{2}}\]
for $r$ large enough.
Further  \[\frac{1}{r+2}{2r+2 \choose r+1}\leqslant \frac{1}{r+1}{2r \choose r}\times 4.\]
We choose $r$ to be the smallest positive integer such that
$n(r+1)\leqslant {2r\choose r}$. We have
\begin{equation}\label{eq:choo}
  n(r+1)\leqslant {2r\choose r}\leqslant 4n(r+1) \text{ and } r\leqslant 3\log n \end{equation}
for $n$ large enough.
We deduce that $\ell = {2r\choose r}-n\leqslant 4n(r+1)\leqslant \cO n\log n$.
So \[\ell^{\ell/2n}\leqslant n^{\cO \log n}.\]
From Equation (\ref{eq:choo}) we deduce
that ${2r\choose r}\leqslant \cO n\log n$.
Also $n^2d(r+1)\leqslant \cO n^2\log^2 n$ and
\[\left(n^2d(r+1)\right)^{r}\leqslant n^{\cO \log n}.\] So   the coefficients of equations
$E_i$ are bounded in absolute value by \[n^{\cO \log n} d_\KK^{\frac{\cO \log n}{n}}.\]
This proves  Theorem \ref{th:sm}.
Theorem \ref{th:nf} follows because there are $r{2r\choose r}$ coefficients to be fixed. We note also that
there may appear several number fields in the smooth zero dimensional part of the complete
intersection $E_1=E_2=\dots=E_r=0$. However the Chow class of this intersection is $r^r\leqslant (\log n)^{\cO\log n}$
and  the number of isolated points  is bounded by this intersection number \cite[Chapter 13]{Fulton}.

\bibliographystyle{plain}
\bibliography{nf}

\end{document}